\documentclass[12pt]{article}

\usepackage{amsfonts}
\usepackage{amsmath}
\usepackage{amssymb}
\usepackage{amscd}
\usepackage{amsthm}
\usepackage{indentfirst}

\newtheorem{thm}{Theorem}[section]
\newtheorem{lemma}[thm]{Lemma}
\newtheorem{prop}[thm]{Proposition}

\theoremstyle{definition}

\newtheorem{rem}[thm]{Remark}
\newtheorem{exam}[thm]{Example}

\newtheorem*{acknow}{Acknowledgements}
\newtheorem*{prf}{Proof}

\newcommand{\R}{{\mathbb{R}}}

\newcommand{\Z}{{\mathbb{Z}}}

\newcommand{\cA}{{\mathcal{A}}}

\newcommand{\cC}{{\mathcal{C}}}

\newcommand{\cM}{{\mathcal{M}}}

\newcommand{\cO}{{\mathcal{O}}}

\newcommand{\cQ}{{\mathcal{Q}}}

\newcommand{\cS}{{\mathcal{S}}}

\newcommand{\fc}{{:\ }}
\newcommand{\ve}{\varepsilon}

\newcommand{\tb}{\textbf}
\newcommand{\ts}{\textsl}

\newcommand{\1}{\ensuremath 1\hspace{-0.25em}{\text{l}}}

\DeclareMathOperator{\im}{im}
\DeclareMathOperator{\id}{id}
\DeclareMathOperator{\area}{area}

\title{Quasi-states and the Poisson bracket on surfaces}
\author{Frol Zapolsky}

\begin{document}

\maketitle

\begin{abstract}

We present a convexity-type result concerning simple quasi-states on closed manifolds. As a corollary, an inequality emerges which relates the Poisson bracket and the measure of non-additivity of a simple quasi-state on a closed surface equipped with an area form. In addition, we prove that the uniform norm of the Poisson bracket of two functions on a surface is stable from below under \(C^0\)-perturbations.

\end{abstract}

\section{Introduction and results}

Let \(X\) be a compact Hausdorff space. We write \(C(X)\) for the Banach algebra of all real-valued continuous functions on \(X\), taken with the supremum norm \(\|\cdot\|\). For \(F \in C(X)\) we denote by \(C(F)\) the closed subalgebra of \(C(X)\) generated by \(F\) and the constant function \(1\), that is
%\[C(F) = \overline{\{p \circ F\,|\, p \in \R[t]\}} = \{\varphi \circ F\,|\, \varphi \in C(\im F)\},\]
\[C(F) = \{\varphi \circ F\,|\, \varphi \in C(\im F)\}.\]
%where \(\R[t]\) is the algebra of polynomials over \(\R\), and the second equality is due to the Weierstrass theorem.

A \textsl{quasi-state} \(\zeta\) on \(X\) is a functional \(\zeta \fc C(X) \to \R\) which satisfies:

(i) \(\zeta(F) \geq 0\) for \(F \geq 0\);

(ii) \(\zeta\) is linear on any \(C(F)\);

(iii) \(\zeta(1) = 1\).\\
A quasi-state \(\zeta\) on \(X\) is \ts{simple} if it is multiplicative on each \(C(F)\).

Quasi-states (as defined here) were introduced and studied by Aarnes, see \cite{quasi-states} and references therein.

Any positive continuous linear functional of norm \(1\) (in other words, a Borel probability measure) furnishes an example of a quasi-state\footnote{Such a quasi-state is simple if and only if it is the \(\delta\)-measure at some point of \(X\).}. However, on many spaces there exist genuine, that is nonlinear, quasi-states (see examples below).

The extent of nonlinearity of a quasi-state \(\zeta\) is measured by the functional \(\Pi(F,G):= \big| \zeta(F+G) - \zeta(F) - \zeta(G) \big|\). It will be clear from the context which quasi-state is meant.

Our first result concerns simple quasi-states on manifolds, namely
\begin{thm}\label{moment_map}
Let \(M\) be a closed manifold, and let \(F,G\) be two continuous functions on \(M\). Let \(\zeta\) be a simple quasi-state on \(M\). Then the image of the moment map \(\Phi \fc M \to \R^2\), \(x \mapsto (F(x), G(x))\), contains the convex hull of the three points \((\zeta(F), \zeta(G))\), \((\zeta(F), \zeta(F+G) - \zeta(F))\), \((\zeta(F+G) - \zeta(G), \zeta(G))\). This is an isosceles right triangle whose legs are of length \(\Pi(F,G)\).
\end{thm}

\begin{rem}
This triangle is in general the largest subset of \(\im \Phi\) we can hope for. For example, if \(M = S^2\) is the round sphere in \(\R^3(x,y,z)\), \(\zeta\) is the (unique) symplectically invariant simple quasi-state on \(M\) (see example \ref{qm_simply_conn} below), and \(F(x,y,z) = x^2\), \(G(x,y,z) = y^2\), then the image of the moment map is the triangle with vertices \((0,0),(0,1),(1,0)\). Here \(\zeta(F) = \zeta(G) = 0\) and \(\Pi(F,G) = \zeta(F+G) - \zeta(F) - \zeta(G) = \zeta(F + G) = 1\).
\end{rem}

\begin{rem}
We do not know whether this result is true for the so-called \ts{pure} quasi-states, that is quasi-states at the extremal boundary of the convex sets of all quasi-states \cite{pure_quasi}.
\end{rem}

This theorem has an interesting corollary in the symplectic context:

\begin{thm}\label{simple_qs_surf}
If \(\zeta\) is a simple quasi-state on a closed surface \(M\) endowed with an area form \(\omega\), then for any \(F,G \in C^\infty(M)\) we have
\begin{equation}\label{simple_qs_surf_ineq}\Pi(F,G)^2 \leq \area(M)\|\{F,G\}\|.\end{equation}
Here \(\area(M) = \int_M \omega\), \(\{\cdot,\cdot\}\) is the Poisson bracket, and \(\|\cdot\|\) stands for the supremum norm.
\end{thm}

\begin{rem}
The theorem easily extends to all representable quasi-states, that is, quasi-states that are elements of the closed convex hull of the set of simple quasi-states, see \cite{pure_quasi}. Indeed, it is clear that \eqref{simple_qs_surf_ineq} holds if \(\zeta\) is a finite convex combination of simple quasi-states. Moreover, let \(\zeta_\nu \to \zeta\) be a net of quasi-states converging to \(\zeta\), such that every \(\zeta_\nu\) satisfies \eqref{simple_qs_surf_ineq}. The topology on the space of quasi-states is such that \(\zeta_\nu(H) \to \zeta(H)\) for any \(H \in C(M)\). Then, since
\[\big|\zeta_\nu(F+G) - \zeta_\nu(F) - \zeta_\nu(G)\big| \leq \area(M) \|\{F,G\}\|,\]
for every \(\nu\), we obtain the same inequality for the limit \(\zeta\).
\end{rem}

\begin{rem}
Note that if \(F,G\) are only required to be of class \(C^1\), the inequality remains valid, for we can choose sequences of \(C^\infty\) functions, \(F_n\), \(G_n\), which tend to \(F,G\), respectively, in the \(C^1\) norm. The claim follows because both sides are continuous with respect to \(C^1\) topology on \(C^\infty(M)\).

In \cite{qs_sympl} it is proved that if \(\zeta\) is a quasi-state on a surface and \(F,G\) are two \(C^\infty\) functions, then \(\{F,G\} \equiv 0 \Rightarrow \Pi(F,G)=0\). In problem 8.2 (ibid.) the authors ask if it is possible to relax the smoothness assumption on \(F,G\), for example to show that if \(\zeta\) is a quasi-state then \(\zeta(F+G) = \zeta(F)+\zeta(G)\) for any \(F,G \in C^1(M)\) with \(\{F,G\} \equiv 0\). The present considerations show that if \(\zeta\) is a representable quasi-state on a closed surface, then it satisfies this property.
\end{rem}

In order to put this result in the proper context, we need the following definition.

Let \((M,\omega)\) be a closed symplectic manifold. A \ts{symplectic quasi-state} \(\zeta\) on \(M\) is a functional \(\zeta \fc C(M) \to \R\) which satisfies:

(i) \(\zeta(F) \geq 0\) for \(F \geq 0\);

(ii) \(\zeta\) linear on Poisson-commutative subalgebras of \(C^\infty(M)\);

(iii) \(\zeta(1) = 1\).

Note that if we require that a functional \(\zeta \fc C(M) \to \R\) satisfy these properties, it automatically becomes linear on every singly generated subalgebra \(C(F)\) of \(C(M)\), therefore every symplectic quasi-state is in particular a quasi-state, and so the terminology is consistent. In dimension two any quasi-state is symplectic, as is proved in \cite{qs_sympl}, and therefore in this case the two notions coincide. The following result appears in \cite{quasimorphism}:

\begin{thm}\label{ineq_EPZ} On certain closed symplectic manifolds \((M,\omega)\) there exist symplectic quasi-states \(\zeta\) which satisfy the following inequality:
\begin{equation}\label{ineq_EPZ_ineq}\Pi(F,G)^2 \leq K(M,\omega)\|\{F,G\}\|,\end{equation}
for any \(F,G \in C^\infty(M)\). Here \(K(M,\omega)\) is a constant depending only on the symplectic manifold.
\end{thm}

A quasi-state is symplectic if \(\{F,G\} = 0 \Rightarrow \Pi(F,G) = 0\). Theorems \ref{simple_qs_surf} and \ref{ineq_EPZ} assert that in case the quantity \(\Pi(F,G)\) is nonzero, it can still be controlled via the Poisson bracket.

One of the manifolds for which the conclusion of theorem \ref{ineq_EPZ} is valid is the standard symplectic sphere, while the corresponding symplectic quasi-state \(\zeta\) is simple. Theorem \ref{simple_qs_surf} then can be viewed as an extension to theorem \ref{ineq_EPZ}, in that it shows that \eqref{ineq_EPZ_ineq} holds (with an appropriate constant) for any closed surface with an area form and any simple quasi-state on it. Also, in the case of the sphere the symplectic quasi-state can be described in elementary terms. However, its origin, as well as the proof of the inequality, both lie in the Floer theory, and are very indirect. One of the motivations for theorem \ref{simple_qs_surf} was to find an elementary proof for \eqref{ineq_EPZ_ineq}. Such a proof has indeed been found, and so this answers part of the question raised in \cite[Section 5]{quasimorphism}.

Another aspect of these inequalities lies in the fact that a quasi-state is Lipschitz, that is \(|\zeta(F) - \zeta(G)| \leq \|F - G\|\) (see \cite{quasi-states}), and so is \(\Pi\): \(|\Pi(F,G) - \Pi(F',G')| \leq 2\big(\|F - F'\|+\|G - G'\|\big)\). Hence the left-hand sides of the inequalities are stable with respect to \(C^0\)-perturbations, while the Poisson bracket, which contains derivatives in its definition, can go wild as a result of such perturbations. But the inequalities tell us that if \(\zeta\) is not additive on a pair of functions, then arbitrarily small \(C^0\)-perturbations cannot make their Poisson bracket vanish. In fact, more can be said. Let us define the following quantity for a pair of smooth functions \(F,G\) on a symplectic manifold \(M\):
\[\Upsilon(F,G) = \liminf_{\ve \to 0} \big\{\|\{F',G'\}\|\,\big|\, F',G' \in C^\infty(M):\, \|F-F'\|,\|G-G'\| < \ve\big\}.\]
A theorem due to Cardin and Viterbo \cite{MTHJ} states that \(\{F,G\} \neq 0\) if and only if \(\Upsilon(F,G) \neq 0\).
Inequality \eqref{ineq_EPZ_ineq} then provides an explicit lower bound on \(\Upsilon(F,G)\) in terms of \(\zeta\) for certain \((M,\omega)\), see \cite{quasimorphism}. Also, in \cite{quasimorphism} the following question was posed: is it true that \(\Upsilon(F,G) = \|\{F,G\}\|\) for any smooth \(F,G\)? As is shown here, in case the manifold is \ts{two-dimensional}, the answer is affirmative:

\begin{thm}\label{eq_Poisson}
Let \((M,\omega)\) be a two-dimensional symplectic manifold (not necessarily closed). For \(F,G \in C^\infty(M)\) we have
\[\Upsilon(F,G) = \|\{F,G\}\|.\]
\end{thm}

\begin{rem}
When \(\{F,G\}\) is an unbounded function, its ``supremum norm'' \(\|\{F,G\}\| = \infty\). It will be clear from the proof that in this case \(\Upsilon(F,G) = \infty\) as well.
\end{rem}

Actually, in the two-dimensional case the Poisson bracket is ``locally stable from below'', see proposition \ref{local_stability} for the precise statement.

\begin{acknow}
I would like to thank my advisor Prof. Leonid Polte\-ro\-vich for arousing my interest in quasi-states, and for his
suggestions, incisive comments and advice. Thanks also to Judy Kupferman for constant supervision of my English, and
for general support and interest. I would like to thank Egor Shelukhin for a suggestion which enabled me to simplify
the proof of lemma \ref{surj}. And finally, I wish to thank Rami Aizenbud for his ongoing curiosity and encouragement
of my work on this topic.\end{acknow}

\section{Definitions and examples}

\subsection{The Poisson bracket}\label{Poisson_br}

We shall employ the following sign convention in the definition of the Poisson bracket of \(F,G \in C^\infty(M)\), where \((M,\omega)\) is a symplectic manifold of dimension \(2n\):
\[-dF \wedge dG \wedge \omega^{n-1} = \textstyle \frac 1 n \{F,G\} \omega^n.\]
%Here \(\omega^k = \omega \wedge \dots \wedge \omega\) is the \(k\)-th exterior power of the symplectic form.
In what follows we shall mostly use the case \(n=1\), and then the formula simplifies to
\[-dF \wedge dG = \{F,G\} \omega.\]
This can be rewritten as follows. If \(\Phi \fc M \to \R^2(x,y)\) is defined by \(\Phi(z) = (G(z),F(z))\) and \(\omega_0 = dx \wedge dy\), then \(\Phi^*\omega_0 = dG \wedge dF = \{F,G\}\omega\).

\subsection{Quasi-states and quasi-measures}

A \ts{space} will always refer to a compact Hausdorff space, unless otherwise mentioned.

We have already defined quasi-states. Let us write \(\cQ(X)\) for the collection of quasi-states on \(X\).

We now turn to another type of objects, called quasi-measures. These are related in a special way to quasi-states, and play a significant role in the present document.

Let \(X\) be a space and let \(\cC\) and \(\cO\) be the collections of closed and open sets in \(X\), respectively. Let
\(\cA = \cC \cup \cO\). A \ts{quasi-measure} \(\tau\) on \(X\) is a function \(\tau \fc \cA \to [0,1]\) which
satisfies:

(i) If \(\{A_i\}_i \subset \cA\) is a finite collection of pairwise disjoint of subsets of \(X\) such that
\(\biguplus_i A_i \in \cA\), then \(\tau \big(\biguplus_i A_i \big) = \sum_i \tau(A_i)\).

(ii) \(\tau(X) = 1\);

(iii) \(\tau(A) \leq \tau(B)\) for \(A,B \in \cA\) such that \(A \subset B\);

(iv) \(\tau(U) = \sup\{\tau(K) \, | \, K \in \cC:\, K \subset U\}\) for \(U \in \cO\).\\
Write \(\cM(X)\) for the collection of quasi-measures on \(X\). A quasi-measure is \ts{simple} if it only takes values
\(0\) and \(1\).

To each quasi-state \(\zeta \in \cQ(X)\) there corresponds a unique quasi-measure \(\tau \in \cM(X)\), defined by the
following formula:
\[\tau(K) = \inf\{\zeta(F) \, | \, F \in C(X): F \geq \1_K\}\]
for \(K \in \cC\), and \(\tau(U) = 1 - \tau(X - U)\) for \(U \in \cO\). Here \(\1_K\) stands for the indicator function
of the set \(K\).

Conversely, to each quasi-measure \(\tau\) there corresponds a unique quasi-state \(\zeta\), obtained through integration with respect to \(\tau\): if \(F \in C(X)\), then the function \(b_F(x) = \tau(\{F < x\})\) is nondecreasing, and takes values in \([0,1]\). Hence it is Riemann integrable, and
\[\zeta(F) \equiv \int_X F\, d\tau = \max_X F - \int\limits_{\min_X F}^{\max_X F} b_F(x)\, dx\,.\]

The described procedures constitute the Aarnes representation theorem, which sets up a bijection \(\cQ(X) \leftrightarrow \cM(X)\) %between the collections of quasi-measures and quasi-states on
for a given space \(X\). This representation theorem is an extension of the Riesz representation theorem, in the sense
that if \(\tau\) is the restriction to \(\cA\) of a Borel probability measure \(\mu\), then the corresponding
quasi-state is the integral with respect to \(\mu\). We refer the reader to \cite{quasi-states} for details.

%Moreover, this bijection is functorial, as we shall now explain. Let \(X,Y\) be two compact Hausdorff spaces and let \(\varphi \fc X \to Y\) be a continuous map. Then \(\varphi\) induces a continuous homomorphism of algebras \(\varphi^* \fc C(Y) \to C(X)\) by \(\varphi^*(F) = F \circ \varphi\), and a map \(\varphi^* \fc \cA(Y) \to \cA(X)\) by \(\varphi^*(A) = \varphi^{-1}(A)\). Also, \(\varphi\) induces maps \(\varphi_* \fc \cQ(X) \to \cQ(Y)\) and \(\varphi_* \fc \cM(X) \to \cM(Y)\), defined by the formulas \(\varphi_*\zeta = \zeta \circ \varphi^*\) and \(\varphi_*\tau = \tau \circ \varphi^*\). Then
%\[\cR(Y) \circ \varphi_* = \varphi_* \circ \cR(X).\]
%Equivalently stated, if \(F \in C(Y)\) and \(\tau \in \cM(X)\), then
%\[\int_X \varphi^* F\, d\tau = \int_Y F\, d(\varphi_* \tau).\]

Another property of this representation theorem is that simple quasi-states correspond to simple quasi-measures \cite{pure_quasi}.

\subsection{Examples of simple quasi-states}

Since simple quasi-states are in bijection with simple quasi-measures, we shall list here examples of the latter.

\begin{exam}\label{qs_diff_inv}
Let \(X\) be a space which is connected and locally connected, and moreover has Aarnes genus \(g = 0\). This last condition is somewhat technical, and since we shall not use it anywhere in the present document, we refer the reader to \cite{qm_construct} for the definition of \(g\) and further details. It suffices to note that in case \(X\) is a compact CW-complex, it has \(g=0\) whenever \(H^1(X;\Z) = 0\) \cite{extreme_qm}. Let us call a subset of \(X\) \ts{solid} if it is connected and has a connected complement.

Let \(\mu\) be a Borel probability measure on \(X\), which has the property that whenever \(K,K'\) are two closed solid
subsets with \(\mu(K) = \mu(K') = \frac 1 2\), then \(K \cap K' \neq \varnothing\). In this case the collection \(\cS\)
of all closed solid subsets of \(X\) having \(\mu\)-measure \(< \frac 1 2\) is a co-basis in the terminology of
\cite{qm_construct}, and so it defines a unique simple quasi-measure \(\tau\), which satisfies
\[\tau(K) = \left\{\begin{array}{ll}0, & \text{if } \mu(K) < \frac 1 2 \\ 1, & \text{otherwise}\end{array}\right. ,\]
for a closed solid \(K\).

We shall mention two particular cases of this construction.

\begin{exam}\label{qm_simply_conn}
Take a simply connected closed manifold \(M\) with a volume form \(\Omega\) satisfying \(\int_M \Omega = 1\), and let \(\mu\) be the Lebesgue measure defined by \(\Omega\). Then the quasi-measure \(\tau\) constructed as above is \(\text{Diff}\,(M,\Omega)\)-invariant.

In case \(M = S^2\) with the standard area form, the resulting simple quasi-state is precisely the one theorem \ref{ineq_EPZ} speaks about.
\end{exam}

\begin{exam}
Take a space \(X\) as in the above example, and let \(\{z_i\}_{i=1}^{2n+1}\) be an odd number of distinct points on \(X\). Let \(\mu = \frac 1 {2n+1} \sum_i \delta_{z_i}\) be the discrete probability measure uniformly distributed among these points. The corresponding quasi-measure may be and is viewed as a generalization of \(\delta\)-measure. See \cite{extreme_qm}.
\end{exam}

\end{exam}

%\begin{exam}\label{qs_pts}
%This is a sort of variation on the theme. Let \(M\) be a simply connected closed manifold, let \(\{z_i\}_{i=1}^{2n+1}\) be an odd number of distinct points on \(M\). Let \(\mu = \frac 1 {2n+1} \sum_i \delta_{z_i}\) be the discrete probability measure uniformly distributed among these points. Then a simple quasi-measure \(\tau\) on \(M\) can be constructed as follows. Suppose that \(K\) is a closed solid subset of \(M\). Then \(\tau(K) = 1\) if the \(\mu\)-measure of \(K\) is at least \(\frac 1 2\) and \(0\) otherwise.
%\end{exam}

\begin{exam}
Examples of simple quasi-measures on the \(2\)-torus have been constructed by Knudsen \cite{extreme_qm}, \cite{qm_torus}.
\end{exam}

\subsection{The median of a Morse function}\label{median}

Let \(M\) be a closed manifold, \(\zeta \in \cQ(M)\) be a simple quasi-state and \(\tau \in \cM(M)\) be the simple
quasi-measure corresponding to \(\zeta\). Let \(F\) be a generic Morse function on \(M\), that is a Morse function with
distinct critical values. The unique component of a level set of \(F\), whose measure with respect to \(\tau\) is
\(1\), is called the \textsl{median of \(F\) relative to \(\zeta\)}, or briefly the median, and is denoted by \(m_F\).
Usually the quasi-state is fixed and so this notation is unambiguous. The median satisfies: \(\zeta(F) = F(m_F)\).

That the median exists can be seen as follows. It is proved in \cite{pure_quasi} that given a continuous function \(G\), the level set \(l_G = G^{-1}(\zeta(G))\) satisfies \(\tau(l_G) = 1\), and is the unique such level set. Now since a level set of a Morse function on a closed manifold is comprised of a finite number of connected components, the finite additivity of \(\tau\) implies the existence of a unique connected component having \(\tau\)-measure \(1\).

The notations introduced here will be used below. The reason for the different notations (\(l\) and \(m\)) is that a continuous function \(G\) need not have a median, and for such a function the only meaningful object is the level set \(l_G\) having quasi-measure \(1\).

\section{Proofs}

\label{proofs}

%\subsection{Existence of the median}

%Maintain the notations of subsection \ref{median}. The level set \(l_F = F^{-1}(\zeta(F))\) satisfies \(\tau(l_F) = 1\), see \cite{pure_quasi}. Now \(\tau\) is finitely additive, and it only takes values \(0\) and \(1\). This, together with the fact that a level set of a generic Morse function has only a finite number of connected components, immediately implies that one of these components, which is of course the median \(m_F\), has quasi-measure \(1\).

\subsection{Proof of theorem \ref{moment_map}}

Denote the triangle mentioned in the theorem by \(\Delta\). Its area is \(\frac 1 2 \Pi(F,G)^2\).

For simplicity assume that \(\zeta(F) = \zeta(G) = 0\), and that \(a := \zeta(F+G) > 0\), in which case \(\Pi(F,G) = \zeta(F+G) = a\). The general case follows easily from this particular one.

Denote \(\kappa(t) = \zeta(F+tG)\). This is a continuous function, and \(\kappa(0) = 0\), \(\kappa(1) = a\). We claim that if \(c,d > 0\) are two numbers such that \(c+d < a\), then the equation \(c+td = \kappa(t)\) has a solution \(t \in (0,1)\). Indeed, this equation can be rewritten as \(c = \kappa(t) - td\). The function on the right side is continuous, and takes values \(0 < c\) for \(t = 0\) and \(a - d > c\) for \(t = 1\). The intermediate value theorem yields the required existence. %This means that the point \((c,d)\) lies on the line \(\{x+ty = \kappa(t)\}\) for some \(t \in (0,1)\).

Now take \(c,d > 0\) such that \(c+d < a\), that is \((c,d) \in \text{Int}\,\Delta\). Fix \(t \in (0,1)\) as above. We shall show that given \(\ve > 0\) there exists a point \(s \in M\) such that \(\|\Phi(s) - (c,d)\| < \ve/t\). Once this is proved, it follows that the image of \(\Phi\) is dense in the triangle \(\Delta\); but \(M\) is compact, and so is its image under the continuous map \(\Phi\), hence \(\im \Phi\) contains the whole of \(\Delta\).

We use the following notation (see subsection \ref{median}): for a continuous function \(E\) on \(M\) let \(l_E := E^{-1}(\zeta(E))\). It follows that \(l_E\) is a set of quasi-measure \(1\). Let \(\ve > 0\) be so small that \(\ve < td\). Put \(H = F+tG\), and let \(K\) be a generic Morse function satisfying \(\|H - K\| < \ve/2\). Denote by \(m_K\) the median of \(K\), as above. Since any two closed sets of quasi-measure \(1\) must intersect, there are points \(p \in m_K \cap l_F\), \(q \in m_K \cap l_G\), \(r \in m_K \cap l_H\). Note that \(\zeta(H) = \kappa(t)\). We have
\begin{align*}
|F(q) - \kappa(t)| &= |(F(q)+tG(q)) - \kappa(t)| \quad &&\text{since } G(q) = \zeta(G) = 0\\
&= |H(q) - H(r)|                                       &&H(r) = \zeta(H) = \kappa(t)\\
&= |H(q) - K(q) + K(r) - H(r)|                         &&\text{since } q,r \in m_K\\
& \leq |H(q) - K(q)| + |H(r) - K(r)|\\
&< \ve.
\end{align*}
In particular, \(c \in (0, F(q))\), since \(c = \kappa(t) - td < F(q) + \ve - td < F(q)\).

Now the points \((0,G(p))\) and \((F(q),0)\) lie in the set \(\Phi(m_K)\) by construction. But \(m_K\) is connected, hence if we denote by \(\pi\) the projection \(\R^2 \to \R,\, (x,y) \mapsto x\), then \(\pi(\Phi(m_K))\), as a connected subset of the real line containing the points \(0\) and \(F(q)\), must contain the entire segment \([0,F(q)]\).

There is then a point \(s \in m_K\) such that \(F(s) = c\). We obtain:
\begin{align*}
t|G(s) - d| &= |(F(s)+tG(s)) - (c+td)| \quad &&\\
&= |H(s) - H(r)|                                       &&H(r) = \kappa(t) = c+td\\
& \leq |H(s) - K(s)| + |H(r) - K(r)|                   &&r,s \in m_K \\
&< \ve.
\end{align*}
Thus
\[\|\Phi(s) - (c,d)\| = \big\|\big(F(s) - c, G(s) - d\big)\big\| = |G(s) - d| < \frac \ve t\,,\]
as required. The proof is thus completed. \qed

\subsection{Proof of theorem \ref{simple_qs_surf}}

In \cite[theorem 3.2.3]{geom_meas} there is proved the so-called area formula. We shall make use of some corollary of it: let \(M\) and \(N\) be two smooth manifolds of dimension \(n\) with \(M\) compact, let \(\Phi \fc M \to N\) be a smooth map, and let \(\Omega\) be a smooth \(n\)-density on \(N\). Then the function \(n_\Phi(z) = \# \Phi^{-1}(z)\), defined on \(N\), is almost everywhere real-valued, and
\[\int_M \Phi^*\Omega = \int_N n_\Phi\Omega.\]

%This is just the coarea formula: if \(M\) and \(N\) are two manifolds of the same dimension, \(M\) is compact, \(\Omega\) is a volume form on \(N\), and \(\Phi \fc M \to N\) is a smooth map, then the function \(n_\Phi(z) = \# \Phi^{-1}(z)\), defined on \(N\), is almost everywhere real-valued, and
%\[\int_M |\Phi^*\Omega| \geq \int_N n_\Phi|\Omega|.\]

In our case \(M\) is the given surface, \(N = \R^2(x,y)\) with the standard density \(\Omega = dx\,dy\), and \(\Phi \fc M \to \R^2\) is \(\Phi(z) = (F(z),G(z))\). It follows from theorem \ref{moment_map} that the image of \(\Phi\) contains a triangle \(\Delta\) of area \(\int_\Delta \Omega = \frac 1 2 \Pi(F,G)^2\). It is true that \(n_\Phi(z) \geq 2\) for almost every \(z \in \Delta\). Indeed, \(M\) is closed, and \(\Phi\) is not onto. Consequently the degree modulo \(2\) of \(\Phi\) is zero, and therefore any regular value must be of even multiplicity. If a regular value is actually attained by \(\Phi\), then its multiplicity is at least two. Note that \(|\{F,G\}||\omega| = |dF \wedge dG| = \Phi^*\Omega\). Putting all this together, and noting that \(\int_M f\omega = \int_M f|\omega|\) for a continuous \(f\), we obtain finally
\begin{multline*}
\Pi(F,G)^2 = 2\int_{\Delta}\Omega \leq \int_{\R^2}n_\Phi \Omega = \int_M \Phi^*\Omega = \int_M|\{F,G\}|\omega \leq \\ \leq \|\{F,G\}\|\cdot\int_M\omega = \area(M)\|\{F,G\}\|.
\end{multline*}\qed

%\begin{rem}
%We need only \(C^1\) smoothness of \(F,G\) in order for the arguments to go through. Indeed, Sard's theorem remains valid for this class of smoothness, provided the manifolds are of the same dimension, as is the case. The arguments presented here rely only on Sard's theorem.
%\end{rem}

\subsection{Proof of theorem \ref{eq_Poisson}}

We shall need two auxiliary results, which are presented below.

Fix a positive integer \(n\). Denote \(B(r) = \{z \in \R^n \, | \, \|z\| < r\}\) for \(r > 0\), where \(\|z\|\) is the Euclidean length of a vector \(z \in \R^n\).

\begin{lemma} \label{surj} Let \(0 < \delta < r\). Consider \(U = B(r)\). Then if \(\Phi \fc \overline U \to \R^n\) is a continuous map which is a \(\delta\)-perturbation in the \(C^0\) norm of the identity map \(\id_{\overline U}\), meaning that\/ \(\sup_{\|z\| \leq r} \|\Phi(z) - z\| < \delta\), then \(\Phi(U)\) contains the ball \(B(r -\delta)\), and if moreover \(\Phi\) is smooth and \(z \in \im \Phi\) is a regular value, then \(\deg_z \Phi = 1\).
\end{lemma}

\begin{prf}
For \(z \in \R^n\) such that \(\|z\| \in [r,r+\delta]\) define
\[t(z) = \frac{\|z\| - r}{\delta}\,,\quad z_0 = \frac {z}{\|z\|}\,r\,,\quad z_1 = \frac {z}{\|z\|}(r+\delta)\,.\]
Clearly \(t(z) = 0\) for \(\|z\| = r\), \(t(z) = 1\) for \(\|z\| = r + \delta\), and \(t(z) \in [0,1]\) for \(\|z\| \in [r,r+\delta]\). Also, \(\|z_0\| = r\), \(\|z_1\| = r+\delta\), and if \(\|z\| = r\) or \(\|z\| = r+\delta\), then \(z = z_0\) or \(z = z_1\), respectively. Extend the definition of \(\Phi\) to the whole of \(\R^n\) by the formula:
\[\Phi(z) = \left\{\begin{array}{ll}
\Phi(z), & \|z\| \leq r \\
(1-t(z))\Phi(z_0) + t(z)z_1, & r \leq \|z\| \leq r+\delta \\
z, & \|z\| \geq r + \delta \\
\end{array}\right. .\]
This extension is clearly continuous, \(\Phi|_{\R^n - V} = \id\), and \(\|\Phi(z) - z\| < \delta\) for all \(z \in \R^n\). For \(\|z\| \geq r + \delta\) and \(\|z\| \leq r\) this is obvious. For  \(\|z\| \in [r,r+\delta]\) we have \(z = (1-t(z))z_0 + t(z)z_1\), and hence
\begin{align*}
\|\Phi(z) - z\| &= \big\|\big[(1-t(z))\Phi(z_0) + t(z)z_1\big] - \big[(1-t(z))z_0 + t(z)z_1\big] \big\|\\
                &= (1-t(z))\|\Phi(z_0) - z_0\|\\
                &<1\cdot \delta = \delta,
\end{align*}
since \(\|\Phi(z_0) - z_0\| < \delta\) by assumption.

Finally, extend \(\Phi\) to a map \(\Phi \fc S^n \to S^n\) by adding \(\infty\) to \(\R^n\) and setting \(\Phi(\infty) = \infty\). This map is continuous and clearly homotopic to \(\id_{S^n}\). Therefore its degree is \(1\), and in particular it is surjective. We have \[\Phi(S^n - U) \cap B(r-\delta) = \varnothing,\]
since \(\|\Phi(z) - z\| < \delta\) for any \(z \in \R^n\) and \(\Phi(\infty) = \infty\). Hence all the points in \(B(r - \delta)\) must come from points of \(U\).

The last assertion follows from the equality of the degree of a smooth map at a regular value and the degree of the map. \qed
\end{prf}

\begin{prop}\label{local_stability} Let \(V\) be an open neighborhood of \(0 \in \R^2(x,y)\), endowed with the standard area form \(\omega_0 = dx \wedge dy\). Let \(F_0,G_0 \in C^\infty(\overline V)\), and suppose that \(\{F_0,G_0\}(0) = 1\). Then for any \(\ve > 0\) there exists \(\delta > 0\) and an open neighborhood \(U\) of \(0\) such that if \(F,G \in C^\infty(\overline V)\) satisfy \(\|F-F_0\|_{\overline V} < \delta, \linebreak \|G-G_0\|_{\overline V} < \delta\), then there exists \(z \in U\) such that \(\{F,G\}(z) > 1- \ve\).
\end{prop}

\begin{prf} Let \(\Phi_0 \fc V \to \R^2\) be defined by\footnote{This order of coordinates is explained by our sign convention, see subsection \ref{Poisson_br}. With this order the map \(\Phi_0\) is orientation-preserving on a neighborhood of the zero, and the function \(\varphi\) introduced here is indeed positive.} \(\Phi_0(z) = (G_0(z),F_0(z))\). There exists \(r > 0\) and a neighborhood \(W\) of \(0\) such that \(\Phi_0 \fc \overline W \to \overline{B(r)}\) is a diffeomorphism. Moreover, if we define the symplectic form \(\omega\) on \(\overline{B(r)}\) by \(\omega = (\Phi_0^{-1})^*\omega_0\), then \(\Phi_0 \fc (\overline W, \omega_0) \to (\overline{B(r)},\omega)\) is a symplectomorphism. There exists a smooth \ts{positive} function \(\varphi\) such that \(\omega = \varphi\, \omega_0\) throughout \(\overline{B(r)}\). We may assume \(r\) to be so small that \(\varphi < 1 + \ve/2\).

Every differential object on \(\overline W\) can be transferred to \(\overline{B(r)}\) by pushing it forward with \(\Phi_0\). In particular, the functions \(G_0\) and \(F_0\) become the coordinates \(x\) and \(y\), the map \(\Phi_0\) becomes the identity map, and if \(F,G\) are smooth functions satisfying the conditions of the proposition, then the map \(\Phi(z) = (G(z),F(z))\) becomes a \(\delta\)-perturbation of the identity map. Therefore we may apply lemma \ref{surj} and conclude that the image of \(\Phi\) contains \(B(r-\delta)\), and moreover, at a regular value \(z\) we have \(\deg_z \Phi = 1\). Then we can write
\[\int_{B(r-\delta)} \omega_0 \leq \int_{\Phi(B(r))}\omega_0.\]
Now the regular values of \(\Phi\) form an open dense subset of \(\Phi(B(r))\). Let \(z\) be such a regular value. By the so-called stack-of-records theorem, there is a small disk \(Y \ni z\) such that \(\Phi^{-1}(Y)\) falls into a finite number of connected components \(\{Y_i\}\), each carried diffeomorphically by \(\Phi\) onto \(Y\). Then
\[\int_Y \omega_0 = \ve_i \int_{Y_i} \Phi^*\omega_0,\]
where \(\ve_i\) is the sign of the Jacobian of \(\Phi\) on \(Y_i\). Since \(\sum \ve_i = \deg_z\Phi = 1\), this implies
\[\int_Y \omega_0 = \deg_z\Phi \int_Y \omega_0 = \sum_i \ve_i \int_{Y_i} \Phi^*\omega_0 = \int_{\Phi^{-1}(Y)} \Phi^*\omega_0.\]
It then follows that
\[ \int_{\Phi(B(r))}\omega_0= \int_{B(r)} \Phi^*\omega_0 = \int_{B(r)} \{F,G\} \omega \leq \max \{F,G\} \int_{B(r)}\omega.\]
But the last integral is
\[\int_{B(r)}\omega = \int_{B(r)}\varphi\,\omega_0 < (1 + \ve/2)\int_{B(r)}\omega_0 = \pi r^2(1 + \ve/2).\]
Now \(\int_{B(r-\delta)} \omega_0 = \pi (r-\delta)^2\), and hence
\[\max\{F,G\} \geq \frac{\pi (r-\delta)^2}{\pi r^2(1 + \ve/2)} > (1 - \ve/2)(1 - 2\delta/r),\]
and if we choose \(\delta < \ve r/4\), we shall obtain the desired inequality. \qed
\end{prf}

%Scaling one of the functions appropriately, we can conclude that if

Returning to the proof of theorem \ref{eq_Poisson}, let \(F,G \in C^\infty(M)\), \(z \in M\), \(\{F,G\}(z) = a\) and \(\ve > 0\). Rescaling one of the functions appropriately and applying the proposition, we can conclude that there is \(\delta > 0\) such that if \(F',G' \in C^\infty(M)\) with \(\|F - F'\|,\|G-G'\| < \delta\), then there exists a point \(z' \in M\) such that \(|\{F',G'\}(z')| > |a| - \ve\). Therefore
\[\|\{F',G'\}\| > |a| - \ve,\]
whence
\[\Upsilon(F,G) \geq |a| - \ve.\]

But since \(\ve\) is arbitrary, we obtain \(\Upsilon(F,G) \geq |a| = |\{F,G\}(z)|\) for any \(z \in M\), and taking the
supremum over \(M\),
\[\Upsilon(F,G) \geq \|\{F,G\}\|.\]
Since the reverse inequality holds trivially, we have the desired result. \qed

\section{Discussion}

There are several directions in which the presented results could be generalized. First of all, theorem \ref{moment_map} speaks about closed manifolds, but the only thing used in the proof is the fact that any continuous function can be approximated by a function whose level sets have only countably many connected components. The countable additivity of quasi-measures, which is established in \cite{additivity_qm}, allows us to define the median of such a function, and then proceed as above. It would be interesting to find spaces other than manifolds with this property.

The second direction is to try and generalize theorem \ref{simple_qs_surf} to arbitrary closed symplectic manifolds and
symplectic quasi-states on them. The methods presented here fail even in the case of a non-representable quasi-state on
a closed surface. And finally, the question posed in \cite{quasimorphism}, namely whether it is true that
\(\Upsilon(F,G) = \|\{F,G\}\|\), is still open in higher dimensions. Apparently some new methods are needed.

\bigskip

\noindent Frol Zapolsky\\
School of Mathematical Sciences\\
Tel Aviv University\\
Tel Aviv 69978, Israel\\
zapolsky@post.tau.ac.il


\begin{thebibliography}{EPZ}

\bibitem[Aa1]{quasi-states} Aarnes, J. F., \textit{Quasi-states and quasi-measures}, Adv. Math. \tb{86} (1991), no. 1,
41--67.

\bibitem[Aa2]{pure_quasi} Aarnes, J. F., \textit{Pure quasi-states and extremal quasi-measures}, Math. Ann. \tb{295} (1993), no. 4,
575--588.

\bibitem[Aa3]{qm_construct} Aarnes, J. F., \textit{Construction of non-subadditive measures and discretization of Borel measures}, Fund. Math. \tb{147} (1993),
213--237.

\bibitem[CV]{MTHJ} Cardin, F., Viterbo, C., \textit{Commuting Hamiltonians and multi-time Hamilton-Jacobi equations}, preprint,
math.SG/0507418.

\bibitem[EP]{qs_sympl} Entov, M., Polterovich, L., \textit{Quasi-states and symplectic intersections}, Comment. Math. Helv. \tb{81} (2006), no. 1,
75--99.

\bibitem[EPZ]{quasimorphism} Entov, M., Polterovich, L., Zapolsky, F., \textit{Quasi-morphisms and the Poisson bracket}, preprint,
math.SG/0605406.

\bibitem[Fe]{geom_meas} Federer, H., \textit{Geometric measure theory}, Die Grundl. der math. Wiss., vol. 153, Springer-Verlag,
1969.

\bibitem[GL]{additivity_qm} Grubb, D. J., LaBerge, T., \textit{Additivity of quasi-measures}, Proc. Amer. Math. Soc. \tb{126} (1998), no. 10,
3007--3012.

\bibitem[Kn1]{extreme_qm} Knudsen, F. F., \textit{Topology and the construction of extreme quasi-measures}, Adv. Math. \tb{120} (1996), no. 2,
302--321.

\bibitem[Kn2]{qm_torus} Knudsen, F. F., \textit{New topological measures on the torus}, Fund. Math. \tb{185} (2005), no. 3,
287--293.

\end{thebibliography}
\end{document}